\documentclass[a4paper,12pt]{article}
\pdfoutput=1

\usepackage[utf8]{inputenc}
\usepackage[english]{babel}
\usepackage{amsthm,amssymb,amsfonts,amsmath}
\usepackage[authoryear]{natbib}
\usepackage{dsfont}
\usepackage{color,fancybox,graphicx}
\usepackage{url}
\usepackage{psfrag}
\usepackage{color}
\usepackage{pifont}
\usepackage{latexsym,pstricks}

\usepackage{hyperref}
\usepackage{mathrsfs}
\usepackage{pst-tree}
\usepackage{fancyhdr}
\usepackage{etoolbox}
\usepackage{setspace}
\usepackage{enumitem}
\usepackage[ ruled, linesnumbered]{algorithm2e}

\newcommand{\T}{\Theta}

\newcommand{\bX}{{\bf X}}
\newcommand{\bx}{{\bf x}}
\newcommand{\bz}{{\bf z}}

\newcommand{\Mtry}{\mathcal{M}_{\textrm{try}}}
\newcommand{\Pfin}{\mathcal{P}_{\textrm{final}}}

\newcommand{\mtry}{\texttt{mtry}}

\DeclareMathOperator*{\argmin}{arg\,min}
\DeclareMathOperator*{\argmax}{arg\,max}

\bibpunct{(}{)}{;}{a}{,}{,}

\usepackage{hyperref}
\hypersetup{
  colorlinks = true,
  urlcolor = blue,
  linkcolor = blue,
  citecolor = blue,
}

\begin{document}
\begin{center}
{\Large
\textbf{\textsf{A Random Forest Guided Tour}}}
\medskip
\medskip
\end{center}

{\bf G\'erard Biau}\\
{\it Sorbonne Universit\'es, UPMC Univ Paris 06, F-75005, Paris, France\\
\& Institut universitaire de France}\\
\href{mailto:gerard.biau@upmc.fr}{gerard.biau@upmc.fr}
\bigskip
\newline
\noindent{\bf Erwan Scornet }\\
{\it Sorbonne Universit\'es, UPMC Univ Paris 06, F-75005, Paris, France}\\
\href{mailto:erwan.scornet@upmc.fr}{erwan.scornet@upmc.fr}\

\begin{abstract}
\noindent{\rm
The random forest algorithm, proposed by L. Breiman in 2001, has been extremely successful as a general-purpose classification and regression method. The approach, which combines several randomized decision trees and aggregates their predictions by averaging, has shown excellent performance in settings where the number of variables is much larger than the number of observations. Moreover, it is versatile enough to be applied to large-scale problems, is easily adapted to various ad-hoc learning tasks, and returns measures of variable importance. The present article reviews the most recent theoretical and methodological developments for random forests. Emphasis is placed on the mathematical forces driving the algorithm, with special attention given to the selection of parameters, the resampling mechanism, and  variable importance measures. This review is intended to provide non-experts easy access to the main ideas.
\medskip

\medskip

}
\end{abstract}

\section{Introduction}
To take advantage of the sheer size of modern data sets, we now need learning algorithms that scale with the volume of information, while maintaining sufficient statistical efficiency. Random forests, devised by L. Breiman in the early 2000s
\citep[][]{Br01},
are part of the list of the most successful methods currently available to handle data in these cases.  This supervised learning procedure, influenced by the early work of \citet{AmGe96}, \citet{Ho}, and \citet{diet},
operates according to the simple but effective ``divide and conquer'' principle: sample fractions of the data, grow a randomized tree predictor on each small piece,  then paste (aggregate) these predictors together.

What has greatly contributed to the popularity of forests is the fact that they can be applied to a wide range of prediction problems and have few parameters to tune. Aside from being simple to use, the method is generally recognized for its accuracy and its ability to deal with small sample sizes and high-dimensional feature spaces. At the same time, it is easily parallelizable and has therefore the potential to deal with large real-life systems. The corresponding \texttt{R} package \texttt{randomForest} can be freely downloaded on the CRAN website (\href{http://www.r-project.org}{http://www.r-project.org}), while a MapReduce \citep[][]{JeSa08} open source implementation called {\it Partial Decision Forests} is available on the Apache Mahout website at \href{https://mahout.apache.org}{https://mahout.apache.org}. This allows the building of forests using large data sets as long as each partition can be loaded into memory.

The random forest methodology has been successfully involved in various practical problems, including a data science hackathon on air quality prediction (\href{http://www.kaggle.com/c/dsg-hackathon}{http://www.kaggle.com/c/dsg-hackathon}), che\-mo\-in\-for\-ma\-ti\-cs \citep[][]{SvLiToCuShFe03}, ecology \citep[][]{PrIvLi06,CuEdBe07}, 3D object recognition \citep[][]{ShFiCoShFiMoKiBl11}, and bioinformatics \citep[][]{microarray}, just to name a few. J. Howard (Kaggle) and M. Bowles (Biomatica) claim in \citet{HoBo12} that {\it ensembles of decision trees---often known as ``random forests''---have been the most successful general-purpose algorithm in modern times}, while H. Varian, Chief Economist at Google, advocates in \citet{Va14} the use of random forests in econometrics.

On the theoretical side, the story of random forests is less conclusive and, despite their extensive use, little is known about the mathematical properties of the method. The most celebrated theoretical result is that of \citet{Br01}, which offers an upper bound on the generalization error of  forests in terms of correlation and strength of the individual trees. This was followed by a technical note \citep[][]{Br04}, which focuses on a stylized version of the original algorithm \citep[see also][]{Br00a, Br00b}. A critical step was subsequently taken by \citet{LiJe06}, who highlighted an interesting connection between random forests and a particular class of nearest neighbor predictors, further developed by \citet{BiDe10}. In recent years, various theoretical studies have been performed \citep[e.g.,][]{Me06, BiDeLu08,IsKo10,Bi12,Ge12,ZhZeKo12}, analyzing more elaborate models and moving ever closer to the practical situation.  Recent attempts towards narrowing the gap between theory and practice include that of \citet{DeMaFr13}, who prove the consistency of a particular online forest, \citet{Wa14} and \citet{MeHo14a}, who study the asymptotic distribution of forests, and \citet{ScBiVe15}, who show that Breiman's \citeyearpar{Br01} forests are consistent in an additive regression framework.

The difficulty in properly analyzing random forests can be explained by the black-box flavor of the method, which is indeed a subtle combination of different components. Among the forests' essential ingredients, both bagging \citep[][]{Br96} and the Classification And Regression Trees (CART)-split criterion \citep{BrFrOlSt84} play critical roles. Bagging (a contraction of bootstrap-aggregating) is a general aggregation scheme, which generates bootstrap samples from the original data set, constructs a predictor from each sample, and decides by averaging. It is one of the most effective computationally intensive procedures to improve on unstable estimates, especially for large, high-dimensional data sets, where finding a good model in one step is impossible because of the complexity and scale of the problem \citep[][]{BuYu02,KlTaSaJo12,WaHaEf13}. As for the CART-split criterion, it originates from the  influential CART program of \citet{BrFrOlSt84}, and is used in the construction of the individual trees to choose the best cuts perpendicular to the axes. At each node of each tree, the best cut is selected by optimizing the CART-split criterion, based on the so-called {\it Gini impurity} (for classification) or the prediction squared error (for regression).

However, while bagging and the CART-splitting scheme play key roles in the random forest mechanism, both are difficult to analyze with rigorous mathematics, thereby explaining why theoretical studies have so far considered simplified versions of the original procedure. This is often done by simply ignoring the bagging step and/or  replacing the CART-split selection by a more elementary cut protocol. As well as this, in Breiman's \citeyearpar{Br01} forests, each leaf (that is, a terminal node) of  individual trees contains a small number of observations, typically between $1$ and $5$.

The goal of this survey is to embark the reader on a guided tour of random forests. We focus on the theory behind the algorithm, trying to give an overview of major theoretical approaches while discussing their inherent pros and cons. For a more methodological review covering applied aspects of random forests, we refer to the surveys by \citet{CrShKo11} and \citet{BoJaKrKo12}. We start by gently introducing the mathematical context in Section 2 and describe in full detail  Breiman's \citeyearpar{Br01} original algorithm. Section $3$ focuses on the theory for a simplified forest model called \textit{purely random forests}, and emphasizes the connections between forests, nearest neighbor estimates and kernel methods.
Section $4$ provides some elements of theory about  resampling mechanisms, the splitting criterion and the mathematical forces at work in Breiman's approach. Section $5$ is devoted to the theoretical aspects of  associated variable selection procedures. Section $6$ discusses various extensions to random forests including online learning, survival analysis and clustering problems. A short discussion follows in Section $7$.

\section{The random forest estimate}

\subsection{Basic principles}

Let us start with a word of caution. The term ``random forests'' is a bit ambiguous. For some authors, it is but a generic expression for aggregating random decision trees, no matter how the trees are obtained. For others, it refers to Breiman's \citeyearpar{Br01} original algorithm. We essentially adopt the second point of view in the present survey.

As mentioned above, the forest mechanism is versatile enough to deal with both supervised classification and regression tasks. However, to keep things simple, we focus in this introduction on regression analysis, and only briefly survey the classification case. Our objective in this section is to provide a concise but mathematically precise presentation of the algorithm for building a random forest. The general framework is nonparametric regression estimation, in which an input random vector $\bX \in \mathcal{X} \subset \mathds{R}^p$ is observed, and the goal is to predict the square integrable random response $Y\in \mathbb R$ by estimating the regression function $m(\bx)=\mathbb E[Y|\bX=\bx]$. With this aim in mind, we assume we are given a training sample $\mathcal D_n=((\bX_1,Y_1), \hdots, (\bX_n,Y_n))$ of independent random variables distributed as the independent prototype pair $(\bX, Y)$. The goal is to use the data set $\mathcal D_n$ to construct an estimate $m_n: \mathcal{X} \to \mathbb R$ of the function $m$.
In this respect, we say that the regression function estimate $m_n$ is (mean squared error) consistent if $\mathbb E [m_n(\bX)-m(\bX)]^2 \to 0$ as $n \to \infty$ (the expectation is evaluated over $\bX$ and the sample $\mathcal D_n$).

A random forest is a predictor consisting of a collection of $M$ randomized regression trees. For the $j$-th tree in the family, the predicted value at the query point $\bx$ is denoted by $m_n(\bx; \Theta_j,\mathcal D_n)$, where $\Theta_1, \hdots,\Theta_M$ are independent random variables, distributed the same as a generic random variable $\Theta$ and independent of $\mathcal D_n$. In practice, the variable $\Theta$ is used to resample the training set prior to the growing of individual trees and to select the successive directions for splitting---more precise definitions will be given later. In mathematical terms, the $j$-th tree estimate takes the form
\begin{align*}
m_n(\bx; \Theta_j,\mathcal D_n) = \sum_{i \in \mathcal{D}^{\star}_n(\Theta_j)} \frac{\mathds{1}_{\bX_i \in A_n(\bx; \Theta_j, \mathcal{D}_n) }Y_i}{N_n(\bx; \Theta_j,\mathcal{D}_n)},
\end{align*}
where $\mathcal{D}_n^{\star}(\Theta_j)$ is the set of data points selected prior to the tree construction, $A_n(\bx; \Theta_j, \mathcal{D}_n)$ is the cell containing $\bx$, and $N_n(\bx; \Theta_j, \mathcal{D}_n)$ is the number of (preselected) points that fall into $A_n(\bx; \Theta_j, \mathcal{D}_n)$.

At this stage, we note that the trees are combined to form the (finite) forest estimate
\begin{equation}
m_{M,n}(\bx; \Theta_1, \hdots, \Theta_M, \mathcal{D}_n)=\frac{1}{M}\sum_{j=1}^M m_n(\bx; \Theta_j,\mathcal D_n).\label{chapitre0_finite_forest}
\end{equation}
In the \texttt{R} package \texttt{randomForest}, the default value of $M$ (the number of trees in the forest) is $\texttt{ntree}=500$. Since $M$ may be chosen arbitrarily large (limited only by available computing resources), it makes sense, from a modeling point of view, to let $M$ tends to infinity,  and consider instead of  (\ref{chapitre0_finite_forest}) the (infinite) forest estimate
\begin{align*}
m_{\infty, n}(\bx; \mathcal{D}_n)=\mathbb E_{\Theta}\left [m_n(\bx;\Theta,\mathcal D_n)\right].
\end{align*}
In this definition, $\mathbb E_{\Theta}$ denotes the expectation with respect to the random parameter $\Theta$, conditional on $\mathcal D_n$.  In fact, the operation ``$M \to \infty$" is justified by the law of large numbers, which asserts that almost surely, conditional on $\mathcal{D}_n$,
$$ \lim\limits_{M \to \infty} m_{M,n}(\bx; \Theta_1, \hdots, \Theta_M, \mathcal{D}_n) = m_{\infty, n}(\bx; \mathcal{D}_n)$$
(see for instance \citealp[][]{Br01}, and \citealp[][]{Sc15a}, for more information on this limit calculation). In the following, to lighten notation we will simply write $m_{\infty, n}(\bx)$ instead of $m_{\infty, n}(\bx;$ $ \mathcal{D}_n)$.

\subsection{Algorithm}
We now provide some insight on how the individual trees are constructed and how randomness kicks in. In Breiman's \citeyearpar{Br01} original forests, each node of a single tree is associated with a hyperrectangular cell. The root of the tree is $\mathcal{X}$ itself and, at each step of the  tree construction, a node (or equivalently its corresponding cell) is split in two parts. The terminal nodes (or leaves), taken together, form a partition of $\mathcal{X}$.

The algorithm works by growing $M$ different (randomized) trees as follows. Prior to the construction of each tree, $a_n$ observations are drawn at random with (or without) replacement from the original data set. These---and only these---$a_n$ observations (with possible repetitions) are taken into account in the tree building. Then, at each cell of each tree, a split is performed by maximizing the CART-criterion (see below) over $\texttt{mtry}$ directions chosen uniformly at random among the $p$ original ones. (The resulting subset of selected coordinates is called $\Mtry$.) Lastly,  construction of  individual trees is stopped when each cell contains less than \texttt{nodesize} points. For any query point $\bx \in \mathcal X$, each regression tree predicts the average of the $Y_i$ (that were among the $a_n$ points) for which the corresponding $\bX_i$ falls into the cell of $\bx$. We draw attention to the fact that growing the tree and making the final estimation only depends on the $a_n$ preselected data points.
Algorithm 1 describes in full detail how to compute a forest's prediction.
%
\begin{algorithm}[!t]
 \KwIn{Training set $\mathcal{D}_n$, number of trees $M>0$, $a_n  \in \{1, \hdots, n \}$, $\texttt{mtry} \in \{1, \hdots, p \}$, $\texttt{nodesize} \in \{1, \hdots, a_n \}$, and $\bx \in \mathcal{X}$.}

 \KwOut{Prediction of the random forest at $\bx$.}

 \For{$j=1, \hdots, M$}{
 Select $a_n$ points, with (or without) replacement, uniformly in $\mathcal{D}_n$. In the following steps, only these $a_n$ observations are used.\

 Set $\mathcal{P} = (\mathcal{X})$ the list containing the cell  associated with the root of the tree.\

 Set $\Pfin = \emptyset$ an empty list.\

\While{ $\mathcal{P} \neq \emptyset $ }{
Let $A$ be the first element of $\mathcal{P}$.\

	\eIf{$A$ contains less than $\texttt{\emph{nodesize}}$ points or if all $\bX_i \in A$ are equal}{
					Remove the cell $A$ from the list $\mathcal{P}$.

					$\Pfin \leftarrow Concatenate (\Pfin , A)$.

					 		}{
		Select uniformly, without replacement, a subset $\Mtry \subset$ $ \{1,\hdots,$ $p\}$ of cardinality $\mtry$.\

		Select the best split in $A$ by optimizing the CART-split criterion along the coordinates in 					$\mathcal{M}_{\textrm{try}}$ \textit{(see text for details)}.\

		Cut the cell $A$ according to the best split. Call $A_L$ and $A_R$ the two resulting cells. \

		Remove the cell $A$ from the list $\mathcal{P}$.\

		$\mathcal{P} \leftarrow Concatenate( \mathcal{P}, A_L, A_R)$.\
				}
}
 Compute the predicted value $m_n(\bx; \Theta_j, \mathcal{D}_n)$ at $\bx$ equal to the  average of the $Y_i$ falling in the cell of $\bx$ in partition $\Pfin$.
 }
 Compute the random forest estimate $m_{M,n}(\bx; \Theta_1, \hdots, \Theta_M, \mathcal{D}_n)$ at the query point $\bx$ according to  (\ref{chapitre0_finite_forest}).
\caption{Breiman's random forest predicted value at $\bx$.}
\end{algorithm}

Algorithm 1 may seem a bit complicated at first sight, but the underlying ideas are simple. We start by noticing that this algorithm has three important parameters:
\begin{enumerate}
\item $a_n  \in \{1, \hdots, n \}$: the number of sampled data points in each tree;
\item $\mtry \in \{1, \hdots, p \}$: the number of possible directions for splitting at each node of each tree;
\item $\texttt{nodesize} \in \{1, \hdots, a_n \}$: the number of examples in each cell below which the cell is not split.
\end{enumerate}
By default, in the regression mode of the \texttt{R} package \texttt{randomForest}, the parameter $\mtry$ is set to $\lceil p/3 \rceil$ ($\lceil \cdot \rceil$ is the ceiling function), $a_n$ is set to $n$, and \texttt{nodesize} is set to $5$. The role and influence of these three parameters on the accuracy of the method will be thoroughly discussed in the next section.

We still have to describe how the CART-split criterion operates. As for now, we consider for the ease of understanding a tree with no subsampling, which uses the entire and original data set $\mathcal{D}_n$ for its construction. Also, we let $A$ be a generic cell and denote by $N_n(A)$ the number of data points falling in $A$. A cut in $A$ is a pair $(j,z)$, where $j$ is some value (dimension) from $\{1, \hdots, p\}$ and $z$  the position of the cut along the $j$-th coordinate, within the limits of $A$. Let $\mathcal{C}_A$ be the set of all such possible cuts in $A$. Then, with the notation $\bX_i = (\bX_i^{(1)}, \hdots, \bX_i^{(p)} )$, for any $(j,z) \in \mathcal{C}_A$, the CART-split criterion takes the form
\begin{align}
L_{{\mbox{\scriptsize reg}}, n}(j,z) = & \frac{1}{N_n(A)} \sum_{i=1}^n (Y_i - \bar{Y}_{A})^2\mathds{1}_{\bX_i \in A} \nonumber \\
&  - \frac{1}{N_n(A)} \sum_{i=1}^n (Y_i - \bar{Y}_{A_{L}} \mathds{1}_{\bX_i^{(j)}  < z} - \bar{Y}_{A_{R}} \mathds{1}_{\bX_i^{(j)} \geq z})^2 \mathds{1}_{\bX_i \in A},  \label{chapitre0_definition_empirical_CART_criterion}
\end{align}
where $A_L = \{ \bx \in A: \bx^{(j)} < z\}$, $A_R = \{ \bx \in A: \bx^{(j)} \geq z\}$, and $\bar{Y}_{A}$ (resp., $\bar{Y}_{A_{L}}$, $\bar{Y}_{A_{R}}$) is the average of the $Y_i$ belonging to $A$ (resp., $A_{L}$, $A_{R}$), with the convention that the average is equal to $0$ when no point $\bX_i$ belongs to $A$ (resp., $A_{L}$, $A_{R}$). For each cell $A$, the best cut $(j_n^{\star},z_n^{\star})$ is selected by maximizing $L_{{\mbox{\scriptsize reg}}, n}(j,z)$ over $\Mtry$ and $\mathcal{C}_A$; that is,
\begin{align*}
(j_n^{\star},z_n^{\star}) \in \argmax\limits_{\substack{j \in \Mtry\\(j,z) \in \mathcal{C}_A }} L_{{\mbox{\scriptsize reg}}, n}(j,z).
\end{align*}
(To remove some of the ties in the argmax, the best cut is always performed in the middle of two consecutive data points.) Let us finally notice that the above optimization program extends effortlessly to the resampling case, by optimizing over the $a_n$ preselected observations instead of the original data set $\mathcal{D}_n$.

Thus, at each cell of each tree, the algorithm chooses uniformly at random $\mtry$ coordinates in $\{1, \hdots, p\}$, evaluates criterion (\ref{chapitre0_definition_empirical_CART_criterion}) over all possible cuts along the directions in $\Mtry$, and returns the best one. The quality measure (\ref{chapitre0_definition_empirical_CART_criterion}) is the criterion used in the influential CART algorithm of \citet{BrFrOlSt84}. This criterion computes the (renormalized) difference between the empirical variance in the node before and after a cut is performed. There are three essential differences between CART  and a tree of Breiman's \citeyearpar{Br01} forest. First of all, in Breiman's forests, the criterion (\ref{chapitre0_definition_empirical_CART_criterion}) is evaluated over a subset  $\Mtry$ of randomly selected coordinates, and not over the whole range $\{1, \hdots, p\}$. Besides, the individual trees are not pruned, and the final cells do not contain more than $\texttt{nodesize}$ observations (unless all data points in the cell have the same $\bX_i$). At last, each tree is constructed on a subset of $a_n$ examples picked within the initial sample, not on the whole sample $\mathcal D_n$; and only these $a_n$ observations are used to calculate the estimation. When $a_n=n$ (and the resampling is done with replacement), the algorithm runs in bootstrap mode, whereas $a_n<n$ corresponds to subsampling (with or without replacement).

\subsection{Supervised classification}
For simplicity, we only consider here the binary classification problem, keeping in mind that random forests are intrinsically capable of dealing with multi-class problems \citep[see, e.g.,][]{microarray}. In this setting \citep{DeGyLu96}, the random response $Y$ takes values in $\{0, 1\}$ and, given $\bX$, one has to guess the value of $Y$. A classifier, or classification rule, $m_n$ is a Borel measurable function of $\bX$ and $\mathcal D_n$ that attempts to estimate the label $Y$ from $\bX$ and $\mathcal D_n$. In this framework, one says that the classifier $m_n$ is consistent if its conditional probability of error
$$L(m_n)=\mathbb P[m_n(\bX)\neq Y] \underset{n\to \infty}{\to} L^{\star},$$
where $L^{\star}$ is the error of the optimal---but unknown---Bayes classifier:
$$
m^{\star} (\bx) = \left\{
\begin{array}{ll}
1 & \mbox{ if  $\mathbb P[Y=1| \bX=\bx] > \mathbb P[Y=0| \bX=\bx]$}\\
0 & \mbox{ otherwise.}
\end{array}
\right.
$$
In the classification context, the random forest classifier is obtained via a majority vote among the classification trees, that is,
\begin{equation*}
m_{M,n}(\bx; \Theta_1, \hdots, \Theta_M, \mathcal{D}_n) = \left\{
\begin{array}{ll}
1 & \mbox{ if $\frac{1}{M}\sum_{j=1}^M m_n(\bx; \Theta_j,\mathcal D_n) > 1/2$}\\
0 & \mbox{ otherwise.}
\end{array}
\right.
\end{equation*}
If a leaf represents region $A$, then a randomized tree classifier takes the simple form
\begin{eqnarray*}
m_n (\bx;\Theta_j,\mathcal D_n) = \left\{
\begin{array}{ll}
1 & \mbox{ if  $\sum_{i \in \mathcal{D}_n^{\star}(\Theta)} \mathbf 1_{\bX_i \in A, Y_i =1} > \sum_{i \in \mathcal{D}_n^{\star}(\Theta)} \mathbf 1_{\bX_i \in A, Y_i =0}$, $\bx\in A$}\\
0 & \mbox{ otherwise,}
\end{array}
\right.
\end{eqnarray*}
where $\mathcal{D}^{\star}_n(\Theta)$ contains the data points selected in the resampling step.
That is, in each leaf, a majority vote is taken over all $(\bX_i,Y_i)$ for which $\bX_i$ is in the same region. Ties are broken, by convention, in favor of class 0. Algorithm 1 can be easily adapted to do classification by modifying the CART-split criterion for the binary setting. To see this, let us consider a single tree with no subsampling step. For any generic cell $A$, let $p_{0,n}(A)$ (resp., $p_{1,n}(A)$) be the empirical probability of a data point in the cell $A$ having label $0$ (resp., label $1$). Then, for any $(j,z) \in \mathcal{C}_A$, the classification CART-split criterion reads
\begin{align*}
L_{{\mbox{\scriptsize class}}, n}(j,z) & = p_{0,n}(A)p_{1,n}(A) - \frac{N_n(A_L)}{N_n(A)} \times p_{0,n}(A_L)p_{1,n}(A_L)\nonumber\\
& \qquad - \frac{N_n(A_R)}{N_n(A)} \times p_{0,n}(A_R)p_{1,n}(A_R).
\end{align*}
This criterion is based on the so-called {\it Gini impurity measure} $2p_{0,n}(A)p_{1,n}(A)$ \citep[][]{BrFrOlSt84}, which has the following simple interpretation. To classify a data point that falls in cell $A$, one uses the rule that assigns a point, uniformly selected from $\{\bX_i \in A: (\bX_i, Y_i)\in \mathcal{D}_n\}$, to label $\ell$ with probability $p_{\ell,n}(A)$, for $j \in \{0,1\}$. The estimated probability that the item has actually label $\ell$ is $p_{\ell,n}(A)$. Therefore the estimated error under this rule is the Gini index $2p_{0,n}(A)p_{1,n}(A)$.

We note that whenever $Y \in \{0,1\}$, optimizing the classification criterion $L_{{\mbox{\scriptsize class}}, n}$ is equivalent to optimizing its regression counterpart $L_{{\mbox{\scriptsize reg}}, n}$. Thus, in this case, the trees obtained with $L_{{\mbox{\scriptsize class}}, n}$ and $L_{{\mbox{\scriptsize reg}}, n}$ are identical. However, the prediction strategy is different: in the classification regime, each tree uses a local majority vote, whereas in regression the prediction is achieved by a local averaging.

When dealing with classification problems, it is usually recommended to set $\texttt{nodesize} = 1$ and $\mtry=\sqrt{p}$ \citep[see, e.g.,][]{LiWi02}.

We draw attention to the fact that regression estimation may also have an interest in the context of dichotomous and multicategory outcome variables (in this case, it is often termed \textit{probability estimation}). For example, estimating outcome probabilities for individuals is important in many areas of medicine, with applications to surgery, oncology, internal medicine, pathology, pediatrics, and human genetics. We refer the interested reader to \citet{MaKrDaMaZi12} and to the survey papers by \citet{KrLiBiKoKoMaZi14} and \citet{KrLiDiHoWeKoZi14}.

\subsection{Parameter tuning}
Literature focusing on  tuning the parameters $M$, \texttt{mtry}, \texttt{nodesize} and $a_n$ is unfortunately rare, with the notable exception of \citet{microarray}, \citet{BeHeAd08}, and \citet{genuer}.
According to \citet{ScKoZi10}, tuning the forest parameters may result in a computational burden, in particular for big data sets, with hundreds and thousands of samples and variables. To circumvent this issue, \citet{ScKoZi10} implement a fast version of the original algorithm, which they name {\it Random Jungle}.

It is easy to see that the forest's variance  decreases as $M$ grows. Thus, more accurate predictions are likely to be obtained by choosing a large number of trees. Interestingly, picking a large $M$ does not lead to overfitting. In effect, following an argument of \citet{Br01}, we have
\begin{align*}
\lim_{n \to \infty}\mathds{E}[m_{M,n}(\bX; \Theta_1, \hdots, \Theta_M) - m(\bX)]^2
 = \mathds{E}[m_{\infty,n}(\bX) - m(\bX)]^2.
\end{align*}
However, the computational cost for inducing a forest increases linearly with $M$, so  a good choice results from a trade-off between computational complexity ($M$ should not be too large for the computations to finish in a reasonable time) and  accuracy ($M$ must be large enough for  predictions to be stable). In this respect, \citet{microarray} argue that the value of $M$ is irrelevant (provided that $M$ is large enough) in a prediction problem involving microarray data sets, where the aim is to classify patients according to their genetic profiles (typically, less than one hundred patients for several thousand genes). For more details we refer the reader to \citet{genuer}, who offer a thorough discussion on the choice of this parameter in various regression problems.
Another  interesting and related approach is by \citet{LaDeDe01}, who propose a simple procedure that determines a priori a minimum number of tree estimates to combine in order to obtain a prediction accuracy level similar to that of a larger forest. Their experimental results show that it is possible to  significantly limit the number of trees.
\medskip

In the \texttt{R} package \texttt{randomForest}, the default value of the parameter \texttt{nodesize} is 1 for classification and 5 for regression. These values are often reported to be good choices \citep[e.g.,][]{microarray}, despite the fact that this is not supported by solid theory. A simple algorithm to tune the parameter \texttt{nodesize} in the classification setting is discussed in \citet{KrScArZi13}.

The effect of $\mtry$ is thoroughly investigated in \citet{microarray}, who show that this parameter has a little impact on the performance of the method, though larger values may be associated with a reduction in the predictive performance. On the other hand, \citet{genuer} claim that the default value of $\mtry$ is either optimal or too small. Therefore, a conservative approach is to take $\mtry$ as large as possible (limited by available computing resources) and  set $\mtry=p$ (recall that $p$ is the dimension of the $\bX_i$). A data-driven choice of $\mbox{\texttt{mtry}}$ is implemented in the algorithm {\it Forest-RK} of \citet{BeHeAd08}.

Let us finally notice that even if there is no theoretical guarantee to support the default values of the parameters, they are nevertheless easy to tune without requiring an independent validation set. Indeed, the procedure accuracy is estimated internally, during the run, as follows. Since each tree is constructed using a different bootstrap sample from the original data, about one-third of the observations are left out of the bootstrap sample and not used in the construction of the $j$-th tree. In this way, for each tree, a test set---disjoint from the training set---is obtained, and averaging over all these left-out data points and over all trees is known as the {\it out-of-bag} error estimate. Thus, the out-of-bag error, computed on the observations set aside by the resampling prior to the tree building, offers a simple way to adjust the parameters without the need of a validation set. \citep[e.g.,][]{KrScArZi13}.

\section{Simplified models and local averaging estimates}

\subsection{Simplified models}
\label{section_simplified_models}
Despite their widespread use, a gap remains between the theoretical understanding of random forests and their practical performance. This algorithm, which relies on complex data-dependent mechanisms, is difficult to analyze and its basic mathematical properties are still not well understood.

As observed by \citet{DeMaDe13}, this state of affairs has led to  polarization between theoretical and empirical contributions to the literature. Empirically focused
papers describe elaborate extensions to the basic random forest framework but come with no clear guarantees. In contrast, most theoretical papers focus on simplifications or stylized versions of the standard algorithm, where the mathematical analysis is more tractable.

\begin{figure}[!!h]
\centering
\begin{tabular}{ccc}
\includegraphics[scale=0.6]{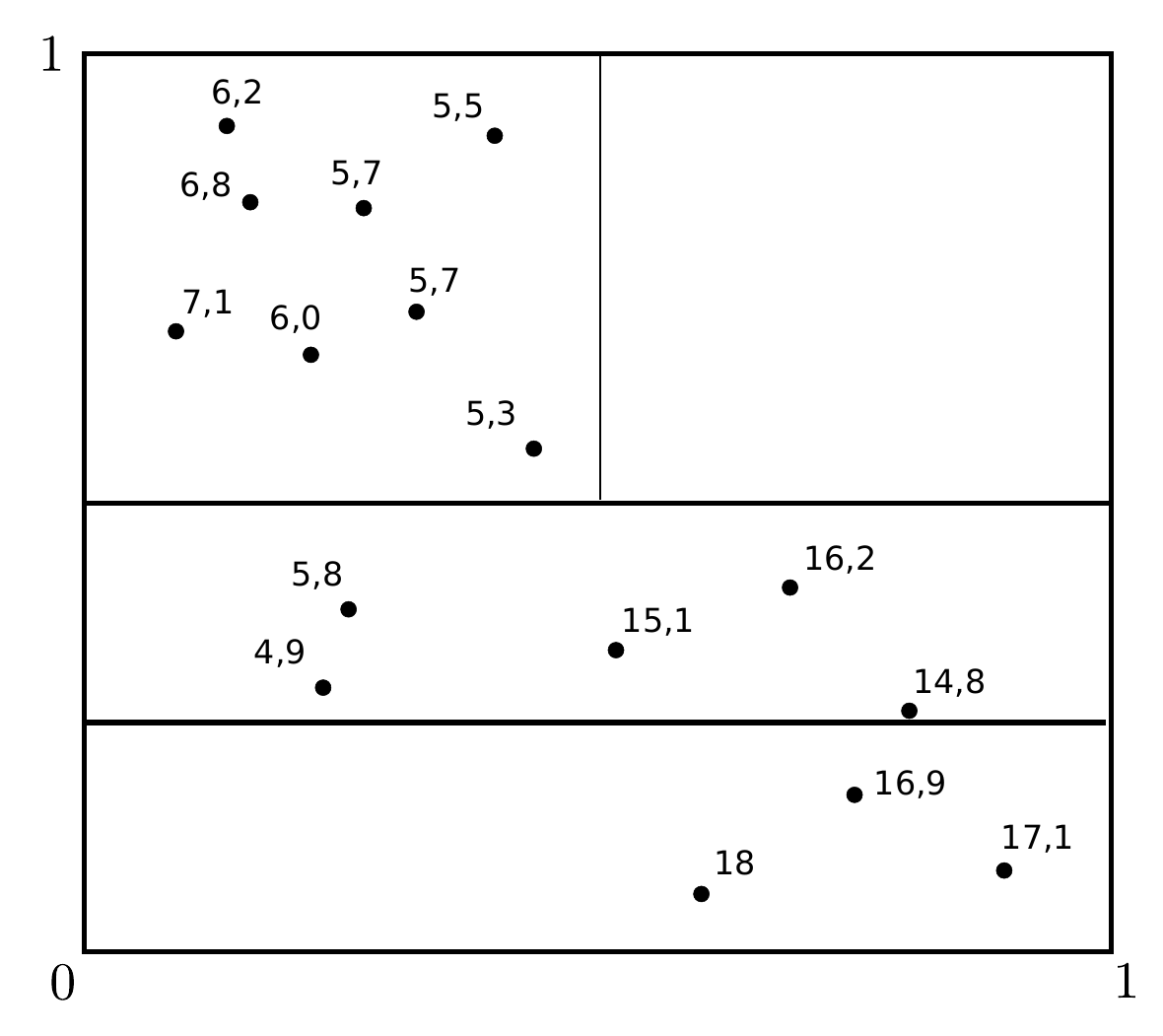}
\end{tabular}
\caption{\label{FIGUREARBRE}A centered tree at level $2$.}
\end{figure}

A basic framework to assess the theoretical properties of forests involves models in which trees are designed independently of the training set $\mathcal{D}_n$. This family of simplified models is often called {\it purely random forests}, for which $\mathcal{X} = [0,1]^d$. A widespread example is the {\it centered forest}, whose principle is as follows: $(i)$ there is no resampling step; $(ii)$ at each node of each individual tree, a coordinate is uniformly chosen in $\{1, \hdots, p\}$; and $(iii)$ a split is performed at the center of the cell along the selected coordinate. The operations $(ii)$-$(iii)$ are recursively repeated $k$ times, where $k\in \mathds{N}$ is a parameter of the algorithm. The procedure stops when a full binary tree with $k$ levels is reached, so that each tree ends up with exactly $2^k$ leaves. The final estimation at the query point $\bx$ is achieved by averaging the $Y_i$ corresponding to the $\bX_i$ in the cell of $\bx$. The parameter $k$ acts as a smoothing parameter that controls the size of the terminal cells (see Figure \ref{FIGUREARBRE} for an example in two dimensions). It should be chosen large enough in order to detect  local changes in the distribution, but not too much to guarantee an effective averaging process in the leaves. In {\it uniform random forests}, a variant of centered forests, cuts are performed uniformly at random over the range of the selected coordinate, not at the center. Modulo some minor modifications, their analysis is similar.

The centered forest rule was first formally analyzed by \citet{Br04}, and then later by \citet{BiDeLu08} and \citet{Sc15a}, who proved that the method is consistent (both for classification and regression) provided $k \to \infty$ and $n/2^k\to \infty$. The proof relies on a general consistency result for random trees stated in \citet[][Chapter 6]{DeGyLu96}. If $\bX$ is uniformly distributed in $\mathcal{X} = [0,1]^p$, then there are on average about $n/2^k$ data points per terminal node. In particular, the choice $k\approx \log n$ corresponds to obtaining a small number of examples in the leaves, in accordance with Breiman's \citeyearpar{Br01} idea that the individual trees should not be pruned. Unfortunately, this choice of $k$ does not satisfy the condition $n/2^k\to \infty$, so  something is lost in the analysis. Moreover, the bagging step is absent, and  forest consistency is obtained as a by-product of  tree consistency. Overall, this model does not demonstrate the benefit of using forests in place of individual trees and is too simple to explain the mathematical forces driving Breiman's forests.

The rates of convergence of centered forests are discussed in \citet{Br04} and \citet{Bi12}. In their approach, the covariates $X^{(j)}$ are independent and the target regression function $m(\bx)=\mathbb E[Y|\bX=\bx]$, which is originally a function of $\bx=(x^{(1)}, \hdots, x^{(p)})$, is assumed to depend only on a nonempty subset $\mathcal{S}$ (for $\mathcal{S}$trong) of the $p$ features. Thus, letting $\bX_{\mathcal{S}}=(X^{(j)}\,:\,j \in \mathcal{S})$,  we have
$$m(\bx)=\mathbb E[Y|\bX_{\mathcal{S}} = \bx_{\mathcal{S}}].$$
The variables of the remaining set $\{1, \hdots, p\} \backslash \mathcal{S}$ have no influence on the function $m$ and can be safely removed. The ambient dimension $p$ can be large, much larger than the sample size $n$, but we believe that the representation is sparse, i.e., that a potentially small number of arguments of $m$ are active--- the ones with indices matching the set $\mathcal{S}$. Letting $|\mathcal{S}|$ be the cardinality of $\mathcal{S}$, the value $|\mathcal{S}|$ characterizes the sparsity of the model: the smaller $|\mathcal{S}|$, the sparser $m$. In this dimension-reduction scenario, \citet{Br04} and \citet{Bi12} proved that if the probability $p_{j,n}$ of splitting along the $j$-th direction tends to $1/S$ and $m$ satisfies a Lipschitz-type smoothness condition, then
\begin{align*}
\mathds{E} \left [m_{\infty, n}(\bX) - m(\bX) \right]^2 =\mbox{O}\Big( n^{\frac{-0.75}{|\mathcal{S}| \log 2 + 0.75}}\Big).
\end{align*}
This equality shows that the rate of convergence of $m_{\infty,n}$ to $m$ depends only on the number $|\mathcal{S}|$ of strong variables, not on the dimension $p$. This rate is strictly faster than the usual rate $n^{-2/(p+2)}$ as soon as $|\mathcal{S}|\leq \lfloor 0.54p\rfloor$ ($\lfloor \cdot \rfloor$ is the floor function). In effect, the intrinsic dimension of the regression problem is $|\mathcal{S}|$, not $p$, and we see that the random forest estimate adapts itself to the sparse framework. Of course, this is achieved by assuming that the procedure succeeds in selecting the informative variables for splitting, which is indeed a strong assumption.

An alternative model for pure forests, called {\it purely uniform random forests} (PURF) is discussed in \citet{Ge12}. For $p=1$, a PURF is obtained by drawing $k$ random variables uniformly on $[0,1]$, and subsequently dividing $[0,1]$ into random sub-intervals. (Note that as such, the PURF can only be defined for $p=1$.). Although this construction is not exactly recursive, it is equivalent to growing a decision tree by deciding at each level which node to split with a probability equal to its length. \citet{Ge12} proves that PURF are consistent and, under a Lipschitz assumption, that the estimate satisfies
\begin{align*}
\mathds{E} [m_{\infty, n}(\bX) - m(\bX) ]^2 =\mbox{O}\left(n^{-2/3}\right).
\end{align*}
This rate is minimax over the class of Lipschitz functions \citep[][]{St80,St82}.

It is often acknowledged that random forests reduce the estimation error of a single tree, while maintaining the same approximation error.
In this respect, \citet{Bi12} argues that the estimation error of centered forests tends to zero (at the slow rate $1/\log n$) even if each tree is fully grown (i.e., $k\approx \log n$). This result is a consequence of the tree-averaging process, since the estimation error of an individual fully grown tree does not tend to zero. Unfortunately, the choice $k\approx \log n$ is too large to ensure consistency of the corresponding forest, whose approximation error  remains constant. Similarly, \citet{Ge12} shows that the estimation error of PURF is reduced by a factor of $0.75$ compared to the estimation error of individual trees. The most recent attempt to assess the gain of forests in terms of estimation and approximation errors is by \citet{ArGe14}, who claim that the rate of the approximation error of certain models is faster than that of the individual trees.
\subsection{Forests, neighbors and kernels}
\begin{figure}[!t]
\psfrag{empty}{empty}
\psfrag{x}{$\bx$}
\centering
\includegraphics*[width=7cm,height=7cm]{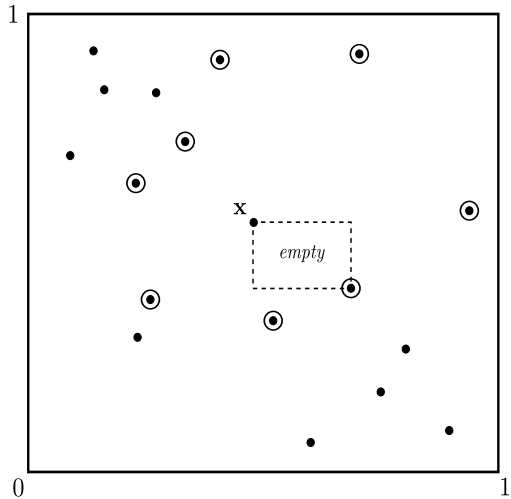}
\caption{\label{figure1} The layered nearest neighbors (LNN) of a point $\bx$ in dimension $p=2$.}
\end{figure}

Let us consider a sequence of independent and identically distributed random variables $\bX_1, \hdots, \bX_n$. In random geometry, an observation $\bX_i$ is said to be a \textit{layered nearest neighbor} (LNN) of a point $\bx$ (from $\bX_1, \hdots, \bX_n$) if the hyperrectangle defined by $\bx$ and $\bX_i$ contains no other data points (\citealp[][]{BaSo66,BaDeHw05}; see also \citealp[][Chapter 11, Problem 6]{DeGyLu96}). As illustrated in Figure \ref{figure1}, the number of LNN of $\bx$ is typically larger than one and depends on the number and configuration of the sample points.

Surprisingly, the LNN concept is intimately connected to  random forests that ignore the resampling step. Indeed, if exactly one point is left in the leaves and if there is no resampling, then no matter what splitting strategy is used, the forest estimate at $\bx$ is a weighted average of the $Y_i$ whose corresponding $\bX_i$ are LNN of $\bx$. In other words,
 \begin{equation}
\label{LNN}
m_{\infty,n}(\bx)=\sum_{i=1}^n W_{ni}(\bx) Y_i,
\end{equation}
where the weights $(W_{n1}, \hdots, W_{nn})$ are nonnegative functions of the sample $\mathcal D_n$ that satisfy  $W_{ni}(\bx)=0$ if $\bX_i$ is not a LNN of $\bx$ and $\sum_{i=1}^n W_{ni}=1$. This important connection was first pointed out by \citet{LiJe06}, who proved that if $\bX$ is uniformly distributed on $[0,1]^p$ then, provided tree growing is independent of $Y_1, \hdots, Y_n$ (such simplified models are sometimes called {\it non-adaptive}), we have
\begin{align*}
\mathds{E} \left[m_{\infty, n}(\bX) - m(\bX)\right ]^2 =\mbox{O}\left( \frac{1}{n_{\mbox{\footnotesize max}} (\log n)^{p-1}}\right),
\end{align*}
where $n_{\mbox{\footnotesize max}}$ is the maximal number of points in the terminal cells (\citealp[][]{BiDe10}, extended this inequality to the case where $\bX$ has a density on $[0,1]^p$). Unfortunately, the exact values of the weight vector $(W_{n1}, \hdots, W_{nn})$ attached to the original random forest algorithm are unknown, and a general theory of forests in the LNN framework is still undeveloped.

It remains however that equation (\ref{LNN}) opens the way to the analysis of random forests via a local averaging approach, i.e., via the average of those $Y_i$ for which $\bX_i$ is ``close'' to $\bx$ \citep[][]{regression}. Indeed, observe starting from (\ref{chapitre0_finite_forest}), that for a finite forest with $M$ trees and without resampling, we have
$$
m_{M,n}(\bx; \T_1, \hdots, \T_M) =  \frac{1}{M} \sum_{j=1}^M \left( \sum_{i=1}^n \frac{Y_i \mathds{1}_{{\bf X}_i \in A_n(\bx; \T_j)}}{N_n(\bx; \Theta_j)}\right),
$$
where $A_n(\bx; \Theta_j)$ is the cell containing $\bx$ and
$
N_n(\bx; \Theta_j) = \sum_{i=1}^n \mathds{1}_{{\bf X}_i \in A_n(\bx; \T_j)}
$
is the number of data points falling in $A_n(\bx; \Theta_j)$. Thus,
$$
m_{M,n}(\bx; \T_1, \hdots, \T_M) =  \sum_{i=1}^n W_{ni}(\bx) Y_i,
$$
where the weights $W_{ni}(\bx)$ are defined by
$$
W_{ni}(\bx) = \frac{1}{M} \sum_{j=1}^M \frac{\mathds{1}_{{\bf X}_i \in A_n(\bx; \T_j)}}{N_n(\bx; \Theta_j)}.
$$
It is easy to see that the $W_{ni}$ are nonnegative and sum to one if the cell containing $\bx$ is not empty. Thus, the contribution of observations falling into cells with a high density of data points is smaller than the contribution of observations belonging to less-populated cells. This remark is especially true when the forests are built independently of the data set---for example, PURF---since, in this case, the number of examples in each cell is not controlled. Next, if we let $M$ tend to infinity, then the estimate $m_{\infty,n}$ may be written (up to some negligible terms)
\begin{equation}
\label{def_kernel_estimates}
m_{\infty,n}(\bx) \approx \frac{\sum_{i=1}^n Y_i K_n(\bX_i, \bx)}{\sum_{j=1}^n K_n(\bX_j, \bx)},
\end{equation}
where
$$K_n(\bx,\bz)=\mathbb P_{\Theta}\left[\bz \in A_n(\bx,\Theta)\right].$$
The function $K_n(\cdot,\cdot)$ is called the {\it kernel} and characterizes the shape of the ``cells'' of the infinite random forest. The quantity $K_n(\bx,\bz)$ is nothing but the probability that $\bx$ and $\bz$ are connected (i.e., they fall in the same cell) in a random tree. Therefore, the kernel $K_n$ can be seen as a proximity measure between two points in the forest. Hence, any forest has its own metric $K_n$, but unfortunately the one associated with the CART-splitting strategy is strongly data-dependent and therefore complicated to work with.

\begin{figure}[h!]
 \hspace{-0.4cm}
 \centering
\begin{tabular}{cc}
\includegraphics[scale=0.33]{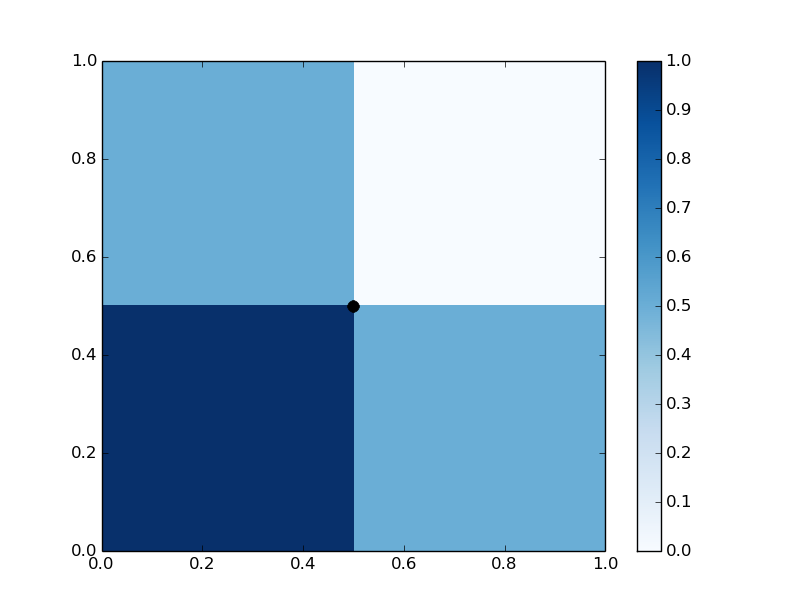}  &
\includegraphics[scale=0.33]{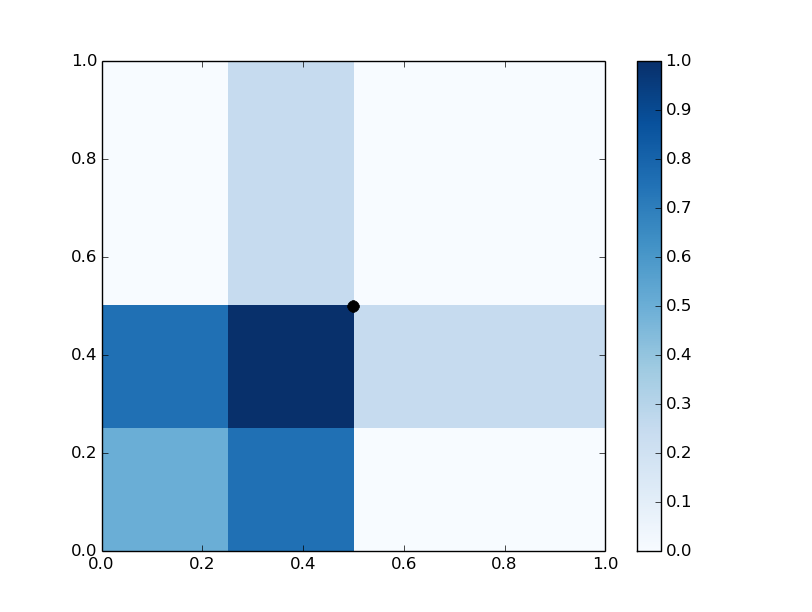}
\end{tabular}
\begin{center}
\includegraphics[scale=0.33]{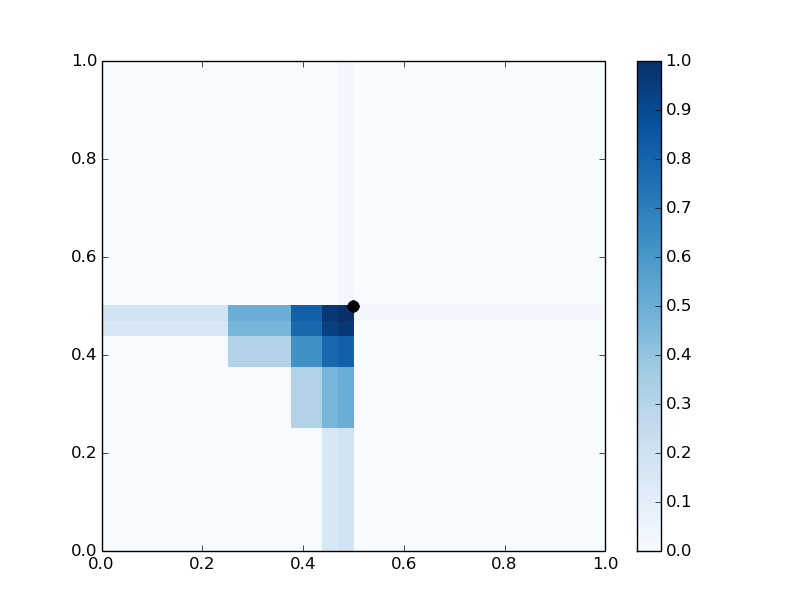}
\end{center}
\caption{Representations of $f_1$, $f_2$ and $f_5$ in $[0,1]^2$. \label{chapitre0_FigureNoyauCentre}}
\end{figure}

It should be noted that $K_n$ does not necessarily belong to the family of Nadaraya-Watson-type kernels \citep[][]{Na64, Wa64}, which satisfy a translation-invariant homogeneous property of the form $K_h(\bx,\bz) = \frac{1}{h}K((\bx - \bz)/h)$ for some {\it smoothing parameter} $h>0$. The analysis of estimates of the form (\ref{def_kernel_estimates}) is, in general, more complicated, depending of the type of forest under investigation. For example, \citet{Sc15b}  proved that for a centered forest defined on $[0,1]^p$ with parameter $k$, we have
$$K_{n,k}(\bx,\bz) = \sum\limits_{\substack{k_{1},\hdots,k_{p} \\ \sum_{j=1}^p k_{j} = k }}
\frac{k!}{k_{1}! \hdots k_{p} !} \left( \frac{1}{p}\right)^k \prod_{j=1}^p  \mathds{1}_{ \lceil 2^{k_j}x_j \rceil = \lceil 2^{k_j}z_j \rceil}.$$
As an illustration, Figure \ref{chapitre0_FigureNoyauCentre} shows the graphical representation for $k=1$, $2$ and $5$ of the function $f_{k}$ defined by
\begin{align*}
\begin{array}{cccl}
f_{k}: 		& 	 [0,1]  \times  [0,1]     & \to    & [0,1]\\
			 & 	 \bz = (   z_1           ,     z_2    ) & \mapsto & K_{n,k}\big((\frac{1}{2},\frac{1}{2}),\bz\big).
\end{array}
\end{align*}

The connection between forests and kernel estimates is mentioned in \citet{Br00a} and developed in detail in \citet{GeErWe06}. The most recent advances in this direction are by \citet{ArGe14}, who show that a simplified forest model can be written as a kernel estimate, and provide its rates of convergence. On the practical side, \citet{DaGh14} plug a specific (random forest-based) kernel---seen as a prior distribution over the piecewise constant functions---into a standard Gaussian process algorithm, and empirically demonstrate that it outperforms the same algorithm ran with linear and radial basis kernels. Besides, forest-based kernels can be used as the input for a large variety of existing kernel-type methods such as Kernel Principal Component Analysis and Support Vector Machines.

\section{Theory for Breiman's forests}

This section deals with Breiman's \citeyearpar{Br01} original algorithm. Since the construction of Breiman's forests depends on the whole sample $\mathcal{D}_n$, a mathematical analysis of the entire algorithm is difficult. To move forward, the individual mechanisms at work in the procedure have been investigated separately, namely the resampling step and the splitting scheme.

\subsection{The resampling mechanism}

The resampling step in Breiman's \citeyearpar{Br01} original algorithm is performed by choosing $n$ times from $n$ points with replacement to compute the individual tree estimates. This procedure, which traces back to the work of \citet[][]{Ef82} \citep[see also][]{PoRoWo99}, is called the {\it bootstrap} in the statistical literature. The idea of generating many bootstrap samples and averaging predictors is called {\it bagging} (bootstrap-aggregating). It was suggested by \cite{Br96} as a simple way to improve the performance of weak or unstable learners. Although one of the great advantages of the bootstrap is its simplicity, the theory turns out to be complex. In effect, the distribution of the bootstrap sample $\mathcal{D}_n^{\star}$ is different from that of the original one $\mathcal{D}_n$, as the following example shows. Assume that $\bX$ has a density, and note that whenever the data points are sampled with replacement then, with positive probability, at least one observation from the original sample is selected more than once. Therefore, with positive probability, there exist two identical data points in $\mathcal{D}_n^{\star}$, and the distribution of $\mathcal{D}_n^{\star}$ cannot be absolutely continuous.
\medskip

The role of the bootstrap in random forests is still poorly understood and, to date, most analyses are doomed to replace the bootstrap by a subsampling scheme, assuming that each tree is grown with $a_n < n$ examples randomly chosen without replacement from the initial sample \citep[][]{MeHo14a,Wa14,ScBiVe15}. Most of the time, the subsampling rate $a_n/n$ is assumed to tend to zero at some prescribed rate---an assumption that excludes de facto the bootstrap regime. In this respect, the analysis of  so-called {\it median random forests} by \citet{Sc15a} provides some insight as to the role and importance of subsampling.

A median forest resembles a centered forest. Once the splitting direction is chosen, the cut is performed at the empirical median of the $\bX_i$ in the cell. In addition, the construction does not stop at level $k$ but continues until each cell contains exactly one observation.
Since the number of cases left in the leaves does not grow with $n$, each tree of a median forest is in general inconsistent \citep[see][Problem $4.3$]{regression}. However, \citet{Sc15a} shows that if $a_n/n \to 0$, then the median forest is consistent, despite the fact that the individual trees are not. The assumption $a_n/n \to 0$ guarantees that every single observation pair $(\bX_i, Y_i)$ is used in the $j$-th tree's construction with a probability that becomes small as $n$ grows. It also forces the query point $\bx$ to be disconnected from  $(\bX_i,Y_i)$ in a large proportion of trees. Indeed, if this were not the case, then the predicted value at $\bx$ would be overly influenced by the single pair $(\bX_i, Y_i)$, which would make the ensemble inconsistent. In fact, the estimation error of the median forest estimate is small as soon as the maximum probability of connection between the query point and all observations is small. Thus, the assumption $a_n/n\to 0$ is but a convenient way to control these probabilities, by ensuring that partitions are dissimilar enough.

\cite{BiDe10} noticed that Breiman's bagging principle has a simple application in the context of nearest neighbor methods.
Recall that the $1$-nearest neighbor (1-NN) regression estimate sets $r_n(\bx)=Y_{(1)}(\bx)$,
where $Y_{(1)}(\bx)$ corresponds to the feature vector $\bX_{(1)}(\bx)$ whose
Euclidean distance to $\bx$ is minimal among all $\bX_1,\ldots,\bX_n$.
(Ties are broken in favor of smallest indices.) It is clearly not, in general, a consistent estimate \citep[Chapter 5]{DeGyLu96}. However, by subbagging, one may turn the $1$-NN estimate into a consistent one, provided that the size of subsamples is sufficiently small. We proceed as follows, via a randomized basic regression estimate $r_{a_n}$ in which $1\leq a_n < n$ is a parameter. The elementary predictor $r_{a_n}$ is the 1-NN rule for a random subsample of size $a_n$ drawn with (or without) replacement from $\mathcal D_n$. We apply subbagging, that is, we repeat the random subsampling an infinite number of times and take the average of the individual outcomes. Thus, the subbagged regression estimate $r^{\star}_n$ is defined by
$$r^{\star}_n(\bx)=\mathds E^{\star}\left[r_{a_n}(\bx)\right],$$
where $\mathbb E^{\star}$ denotes expectation with respect to the resampling distribution, conditional on the data set $\mathcal D_n$.
\cite{BiDe10} proved that the estimate $r^{\star}_n$ is universally  (i.e., without conditions on the distribution of $(\bX,Y)$)  mean squared consistent, provided $a_n\to \infty$ and $a_n/n\to 0$. The proof relies on the observation that $r_n^{\star}$ is in fact a local averaging estimate \citep[][]{St77} with weights
$$W_{ni}(\bx)=\mathbb P[\bX_i \mbox{ is the 1-NN of $\bx$ in a random selection of size } a_n].$$
The connection between bagging and nearest neighbor estimation is further explored by \citet{BiCeGu10}, who prove that the subbagged estimate $r^{\star}_n$ achieves optimal rate of convergence over Lipschitz smoothness classes, independently from the fact that resampling is done with or without replacement.

\subsection{Decision splits}

The coordinate-split process of the random forest algorithm is not easy to grasp, essentially because it uses both the $\bX_i$ and  $Y_i$ variables to make its decision. Building upon the ideas of \citet{BuYu02}, \citet{BaKe07} establish a limit law for the split location in the context of a regression model of the form $Y = m(\bX) + \varepsilon$, where $\bX$ is real-valued and $\varepsilon$  an independent Gaussian noise. In essence, their result is as follows. Assume for now that the distribution of $(\bX,Y)$ is known, and denote by $d^{\star}$ the (optimal) split that maximizes the theoretical CART-criterion at a given node. In this framework, the regression estimate restricted to the left (resp., right) child of the cell takes the form
\begin{align*}
\beta_{\ell,n}^{\star} = \mathds{E} [Y | X \leq d^{\star}] \quad \Big(\textrm{resp.,}~ \beta_{r,n}^{\star} = \mathds{E} [Y | X > d^{\star}] \Big).
\end{align*}
When the distribution of $(\bX,Y)$ is unknown, so are $\beta_{\ell}^{\star}$, $\beta_r^{\star}$ and $d^{\star}$, and these quantities are estimated by their natural empirical counterparts:
\begin{align*}
(\hat{\beta}_{\ell,n}, \hat{\beta}_{r,n}, \hat{d}_n) \in \argmin_{\beta_{\ell}, \beta_r, d}
\sum_{i=1}^n \big[ Y_i - \beta_{\ell} \mathds{1}_{X_i \leq d} - \beta_r \mathds{1}_{X_i > d} \big]^2.
\end{align*}
Assuming that the model satisfies some regularity assumptions (in particular, $\bX$ has a density $f$, and both $f$ and $m$ are continuously differentiable), \citet{BaKe07} prove that
\begin{align}
n^{1/3} \left(
\begin{array}{c}
\hat{\beta}_{\ell,n} - \beta_{\ell}^{\star}\\
\hat{\beta}_{r,n} - \beta_r^{\star} \\
\hat{d}_n - d^{\star}
\end{array}
\right)
\overset{\mathcal{D}}{\to}
\left(
\begin{array}{c}
c_1\\
c_2 \\
1
\end{array}
\right)
\argmax_t (a W(t) - bt^2),\label{limiting_distribution}
\end{align}
where  $\mathcal{D}$ denotes convergence in distribution, and $W$ is a standard two-sided Brownian motion process on the real line. Both $a$ and $b$ are positive constants that depend upon the model parameters and the unknown quantities  $\beta_{\ell}^{\star}$, $\beta_r^{\star}$ and $d^{\star}$. The limiting distribution in (\ref{limiting_distribution}) allows one to construct confidence intervals for the position of CART-splits. Interestingly, \citet{BaKe07} refer to the study of \citet{QiKiRi03} on the effects of phosphorus pollution in the Everglades, which uses split points in a novel way. There, the authors identify threshold levels of phosphorus concentration that are associated with declines in the abundance of certain species. In their approach, split points are not just a means to build trees and forests, but can also provide important information on the structure of the underlying distribution.
\medskip

A further analysis of the behavior of forest splits is performed by \citet{Is13}, who argues that the so-called {\it End-Cut Preference} (ECP) of the CART-splitting procedure (\citealp[that is, the fact that splits along non-informative variables are likely to be near the edges of the cell---see][]{BrFrOlSt84}) can be seen as a desirable property. Given the randomization mechanism at work in  forests, there is indeed a positive probability that none of the preselected variables at a node are informative. When this happens, and if the cut is performed, say, at the center of a side of the cell (assuming that $\mathcal{X} = [0,1]^d$), then the sample size of the two resulting cells is drastically reduced by a factor of two---this is an undesirable property, which may be harmful for the prediction task. Thus, \citet{Is13} stresses that the ECP property ensures that a split along a noisy variable is performed near the edge, thus maximizing the tree node sample size and making it possible for the tree to recover from the split downstream. \citet{Is13} claims that this property can be of benefit even when considering a split on an informative variable, if the corresponding region of space contains little signal.

It is shown in \citet{ScBiVe15} that random forests asymptotically perform, with high probability, splits along the $S$ informative variables (in the sense of Section \ref{section_simplified_models}). Denote by $j_{n,1}(\bX), \hdots,$ $j_{n,k}(\bX)$ the first $k$ cut directions used to construct the cell of $\bX$, with the convention that $j_{n,q}(\bX) = \infty$ if the cell has been cut strictly less than $q$ times. Assuming some regularity conditions on the regression model, and considering a modification of Breiman's forests in which all directions are preselected for splitting, \citet{ScBiVe15} prove that, with probability $1-\xi$, for all $n$ large enough and all $1 \leq q \leq k$,
\begin{align*}
j_{n,q}(\bX) \in \{1, \hdots, S\}.
\end{align*}
This result offers an interesting perspective on why random
forests nicely adapt to the sparsity setting. Indeed, it shows that the algorithm selects splits mostly along the $S$ informative variables, so that everything happens as if data were projected onto the vector space spanned by these variables.

There exists a variety of random forest variants based on the CART-criterion. For example, the {\it Extra-Tree} algorithm of  \citet{GeErWe06} consists in randomly selecting a set of split points and then choosing the split that maximizes the CART-criterion. This algorithm has similar accuracy performance while being more computationally efficient. In the {\it PERT} (Perfect Ensemble Random Trees) approach of \citet{CuZh01}, one builds perfect-fit classification trees with random split selection. While individual trees clearly overfit, the authors claim that the whole procedure is eventually consistent since all classifiers are believed to be almost uncorrelated.
As a variant of the original algorithm, \citet{Br01} considered splitting along linear combinations of features \citep[this procedure has been  implemented by][in the package \texttt{obliquetree} of the statistical computing environment \texttt{R}]{Tr09}.
As noticed by \citet{MeKeSpKoHa11}, the feature space separation by orthogonal hyperplanes in random forests results in box-like decision surfaces, which may be advantageous for some data but suboptimal for other, particularly for collinear data with correlated features .

With respect to the tree building process, selecting uniformly at each cell a set of features for splitting is simple and convenient, but such procedures inevitably select irrelevant variables. Therefore, several authors have proposed modified versions of the algorithm that incorporate a data-driven weighing of  variables. For example, \citet{KyZo14} study the effectiveness of non-uniform randomized feature selection in classification tree, and experimentally show that such an approach may be more effective compared to naive uniform feature selection. {\it Enriched Random Forests}, designed by \citet{AmCaLe08} choose at each node the eligible subsets by weighted random sampling with the weights tilted in favor of  informative features. Similarly, the {\it Reinforcement Learning Trees} (RLT) of \citet{ZhZeKo12} build at each node a random forest to determine the variable that brings the greatest future improvement in later splits, rather than choosing the one with largest marginal effect from the immediate split.
Splits in random forests are known to be biased toward covariates with many possible splits \citep[][]{BrFrOlSt84,Se88} or with missing values \citep[][]{KiLo01}. \citet{HoHoZe06} propose a two-step procedure to correct this situation by first selecting the splitting variable and then the position of the cut along the chosen variable. The predictive performance of the resulting trees is empirically shown to be as good as the performance of the exhaustive search procedure. We also refer the reader to \citet{ZiKo14}, who review the different splitting strategies.

Choosing weights can also be done via regularization. \citet{DeRu12} propose a {\it Regularized Random Forest} (RRF), which penalizes selecting a new feature for splitting when its gain is similar to the features used in previous splits. \citet{DeRu13} suggest a {\it Guided RRF} (GRRF), in which the importance scores from an ordinary random forest are used to guide the feature selection process in RRF. Lastly, a Garrote-style convex penalty, proposed by \citet{Me09}, selects functional groups of nodes in  trees, yielding to parcimonious estimates. We also mention the work of \citet{KoGa14} who address the problem of controlling the false positive rate of random forests and present a principled way to determine thresholds for the selection of relevant features without any additional computational load.

\subsection{Consistency, asymptotic normality, and more}

All in all, little has been proven mathematically for the original procedure of Breiman.
A seminal result by \citet{Br01} shows that the error of the forest is small as soon as the predictive power of each tree is good and the correlation between the tree errors is low. More precisely, independently of the type of forest, one has
\begin{align*}
\mathds{E}_{\bX, Y} [Y - m_{\infty,n}(\bX)]^2
\leq \bar{\rho}\, \mathds{E}_{\Theta,\bX, Y} [Y - m_{n}(\bX; \Theta)]^2,
\end{align*}
where
\begin{align*}
\bar{\rho} = \frac{\mathds{E}_{\Theta,\Theta'}[\rho(\Theta, \Theta') g(\Theta) g(\Theta')]}{\mathds{E}_{\Theta}[g(\Theta)]^2},
\end{align*}
with $\Theta$ and $\Theta'$ independent and identically distributed,
$$
 \rho(\Theta, \Theta') = \textrm{Corr}_{\bX, Y}\big[Y-m_n(\bX; \Theta),Y-m_n(\bX; \Theta')\big],
 $$
 and $g(\Theta) = \sqrt{\mathds{E}_{\bX, Y}[Y - m_n(\bX; \Theta)]^2}$. Similarly, \citet{FrHaTi01} decompose the variance of the forest as a product of the correlation between trees and the variance of a single tree. Thus, for all $\bx$,
\begin{align*}
\textrm{Var}[m_{\infty,n}(\bx)] = \rho (\bx) \sigma (\bx),
\end{align*}
where $\rho (\bx) = \textrm{Corr}[m_n(\bx; \Theta), m_n(\bx; \Theta') ]$ and $\sigma(\bx) = \textrm{Var}[m_n(\bx; \Theta)]$.

A link between the error of the finite and infinite forests is established in \citet{Sc15a}, who shows, provided some regularity assumptions are satisfied, that
\begin{align*}
0 & \leq \mathds{E} [m_{M,n}(\bX; \Theta_1, \hdots, \Theta_M)-m(\bX)]^2 - \mathds{E} [m_{\infty,n}(\bX)-m(\bX)]^2 \\
& \leq \frac{8}{M} \times \big( \|m\|_{\infty}^2 +  \sigma^2 (1 + 4\log n) \big).
\end{align*}
This inequality provides an interesting solution for choosing the number of trees, by making the error of the finite forest arbitrary close to that of the infinite one.

Consistency and  asymptotic normality of the whole algorithm were recently proved, replacing bootstrap by subsampling and simplifying the splitting step. So, \citet{Wa14} shows the asymptotic normality of Breiman's infinite forests, assuming that $(i)$ cuts are spread along all the $p$ directions and do not separate a small fraction of the data set; and $(ii)$ two different data set are used to respectively build the tree and estimate the value within a leaf. He also establishes that the infinitesimal jackknife \citep[][]{Ef79} consistently estimates the forest variance.

\citet{MeHo14a} prove a similar result for finite forests, which mainly relies on the assumption that the prediction of the forests does not much vary when the label of one point in the training set is slightly modified. These authors show that whenever $M_n$ (the number of trees) is allowed to vary with $n$, and when $a_n = \mbox{o}(\sqrt{n})$ and $\lim_{n \to \infty} n/M_n= 0$, then, for a fixed $\bx$,
$$\frac{\sqrt{n}(m_{M,n}(\bx; \Theta_1, \hdots, \Theta_M) - \mathds{E}[m_{\infty,n}(\bx)])}{\sqrt{a_n^2 \zeta_{1,a_n}}} \overset{\mathcal{D}}{\to}N,$$
where $N$ is a standard normal random variable,
$$\zeta_{1,a_n} = \textrm{Cov}\left[m_n(\bX_1, \bX_2, \hdots, \bX_{a_n}; \Theta), m_n(\bX_1, \bX_2', \hdots, \bX_{a_n}'; \Theta')\right],$$
$\bX_i'$  an independent copy of $\bX_i$ and $\Theta'$ an independent copy of $\Theta$. It is worth noting that both \citet{MeHo14a} and \citet{WaHaEf13} provide corrections for estimating the forest variance $\zeta_{1,a_n}$.

\citet{ScBiVe15} proved a consistency result in the context of additive regression models for the pruned version of Breiman's forest. Unfortunately, the consistency of the unpruned procedure comes at the price of  a conjecture regarding the behavior of the CART algorithm that is difficult to verify.

\begin{figure}[t]
\centering
\leavevmode
\includegraphics[scale=0.5]{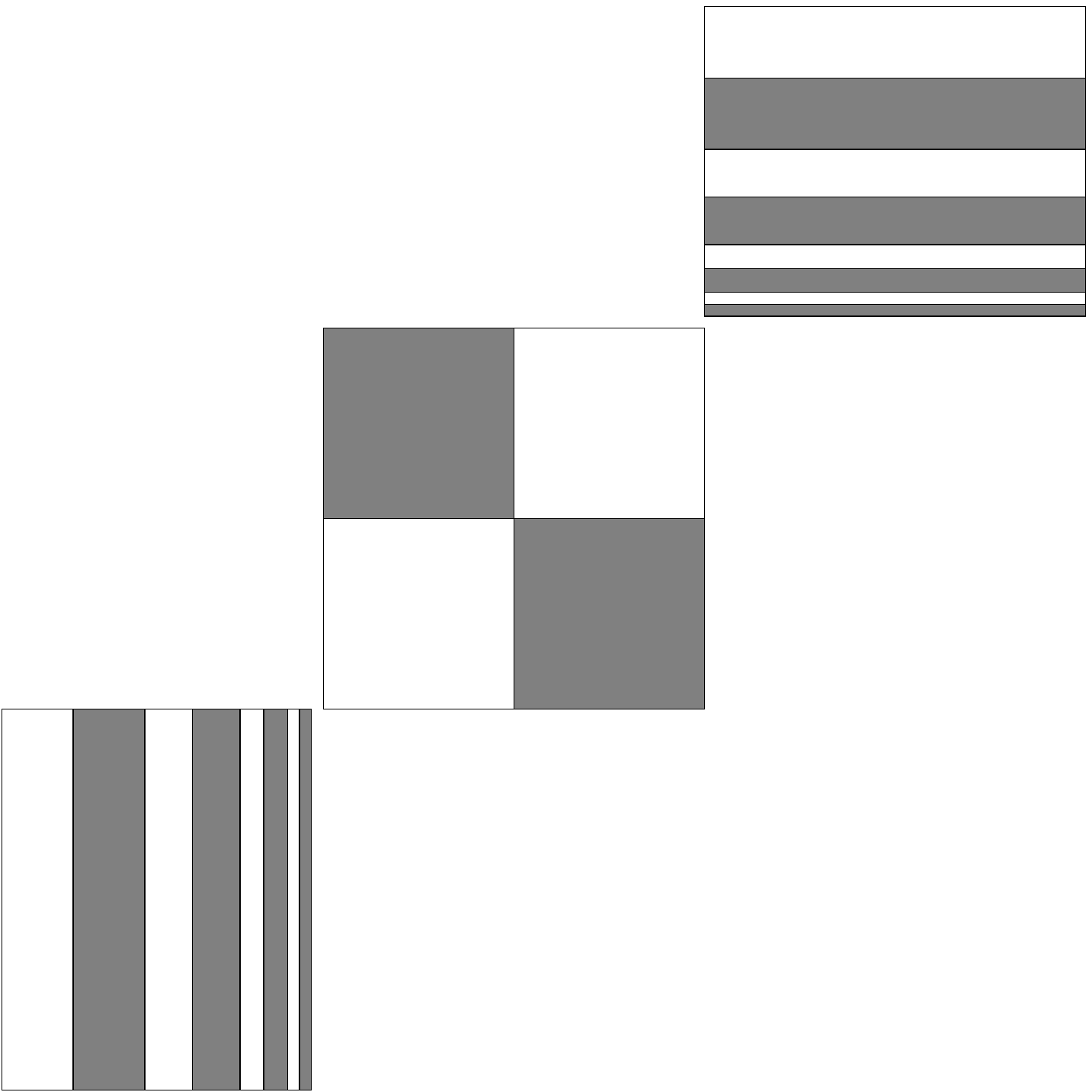}
\caption{
\label{greedy}
An example of a distribution for which greedy random forests
are inconsistent. The distribution of $\bX$ is uniform on the union of
the three large squares. White areas represent the set where $m(\bx)=0$
and grey where $m(\bx)=1$.
}
\end{figure}

We close this section with a negative but interesting result due to \citet{BiDeLu08}. In this example, the total
number $k$ of cuts is fixed and $\texttt{mtry}=1$. Furthermore, each tree is built by minimizing the true probability of error at each node. Consider the joint distribution of $(\bX,Y)$ sketched in Figure \ref{greedy} and let $m(\bx) = \mathds{E}[Y | \bX = \bx]$.
The variable $\bX$ has a uniform distribution on
$[0,1]^2 \cup [1,2]^2 \cup [2,3]^2$ and
$Y$ is a function of $\bX$---that is, $m(\bx) \in \{0,1\}$ and $L^{\star}=0$---defined as follows.
The lower left square $[0,1]\times [0,1]$ is divided into countably
infinitely many vertical strips in which the strips with $m(\bx)=0$
and $m(\bx)=1$ alternate. The upper right square $[2,3]\times [2,3]$
is divided similarly into horizontal strips. The middle rectangle
$[1,2]\times [1,2]$ is a $2\times 2$ checkerboard. It is easy to see that no matter what the sequence of random selection of split directions is and no matter for how long each tree is grown, no tree will ever cut the middle rectangle and therefore the probability of error of the corresponding random forest classifier is at least $1/6$.
This example illustrates that consistency of greedily grown random forests is a delicate issue. We note however that if Breiman's  \citeyearpar{Br01} original algorithm is used in this example (that is, each cell contains exactly one data point) then one obtains a consistent classification rule. We also note that the regression function $m$ is not Lipschitz---a smoothness assumption on which many results on random forests rely.

\section{Variable importance}

\subsection{Variable importance measures}
Random forests can be used to rank the importance of variables in regression or classification problems via two measures of significance.
The first, called {\it Mean Decrease Impurity} \citep[MDI; see][]{Br02}, is based on the total decrease in node impurity from splitting on the variable, averaged over all trees. The second, referred to as {\it Mean Decrease Accuracy} (MDA), first defined by \citet{Br01}, stems from the idea that if the variable is not important, then rearranging its values should not degrade prediction accuracy.

Set $\bX = (X^{(1)}, \hdots, X^{(p)} )$. For a forest resulting from the aggregation of $M$ trees, the MDI of the variable $X^{(j)}$ is defined by
\begin{align*}
\widehat{\mbox{MDI}}(X^{(j)}) = \frac{1}{M} \sum_{{\ell}=1}^M \sum_{\substack{t \in \mathcal{T}_{\ell}\\ j_{n,t}^{\star} = j}}  p_{n,t} L_{{\scriptsize \mbox{reg}}, n}(j_{n,t}^{\star}, z_{n,t}^{\star}),
\end{align*}
where $p_{n,t}$ is the fraction of observations falling in the node $t$, $\{\mathcal{T}_{\ell}\}_{1\leq {\ell} \leq M}$  the collection of trees in the forest, and $(j_{n,t}^{\star}, z_{n,t}^{\star})$ the split that maximizes the empirical criterion (\ref{chapitre0_definition_empirical_CART_criterion}) in  node $t$. Note that the same formula holds for classification random forests by replacing the criterion $L_{{\scriptsize \mbox{reg}}, n}$ by its classification version $L_{{\scriptsize \mbox{class}}, n}$. Thus, the MDI of $X^{(j)}$ computes the weighted decrease of impurity corresponding to splits along the variable $X^{(j)}$ and averages this quantity over all trees.
\medskip

The MDA relies on a different principle and uses the out-of-bag error estimate (see Section $2.4$). To measure the importance of the $j$-th feature, we randomly permute the values of variable $X^{(j)}$ in the out-of-bag observations and put these examples down the tree. The MDA of $X^{(j)}$ is obtained by averaging the difference in out-of-bag error estimation before and after the permutation over all trees. In mathematical terms, consider a variable $X^{(j)}$ and denote by $\mathcal{D}_{{\ell},n}$ the out-of-bag data set of the ${\ell}$-th tree and ${\mathcal{D}}^j_{{\ell},n}$ the same data set where the values of $X^{(j)}$ have been randomly permuted. Recall that $m_n(\cdot; \Theta_{\ell})$ stands for the ${\ell}$-th tree estimate. Then, by definition,
\begin{align}
\widehat{\mbox{MDA}}(X^{(j)})= \frac{1}{M} \sum_{{\ell}=1}^M \bigg[ R_n\big[m_{n}(\cdot; \Theta_{\ell}),{\mathcal{D}}^j_{{\ell},n}\big] - R_n\big[m_{n}(\cdot; \Theta_{\ell}),\mathcal{D}_{{\ell},n}\big] \bigg], \label{test_bis}
\end{align}
where $R_n$ is defined for $\mathcal D=\mathcal{D}_{{\ell},n}$ or $\mathcal D={\mathcal{D}}^j_{{\ell},n}$ by
\begin{align}
R_n\big[m_{n}(\cdot; \Theta_{\ell}), \mathcal{D}\big] = \frac{1}{|\mathcal{D}|} \sum_{i: (\bX_i, Y_i) \in \mathcal{D}} (Y_i - m_{n}(\bX_i; \Theta_{\ell}) )^2.
\label{test_bisb}
\end{align}
It is easy to see that the population version of $\widehat{\mbox{MDA}}(X^{(j)}) $ is
\begin{align*}
\mbox{MDA}^{\star}(X^{(j)}) = \mathds{E} \big[ Y - m_n(\bX'_{j}; \Theta) \big]^2 - \mathds{E} \big[ Y - m_n(\bX;\Theta) \big]^2,
\end{align*}
where $\bX'_{j} = (X^{(1)}, \hdots, X'^{(j)} , \hdots, X^{(p)})$ and $X'^{(j)}$ is an independent copy of $X^{(j)}$. For classification purposes, the MDA still satisfies (\ref{test_bis}) and (\ref{test_bisb}) since $Y_i \in \{0,1\}$ (so, $R_n(m_n(\cdot; \Theta),\mathcal{D})$ is also the proportion of points that are correctly classified by $m_n(\cdot; \Theta)$ in $\mathcal{D}$).

\subsection{Theoretical results}
In the context of a pair of categorical variables $(\bX,Y)$, where $\bX$ takes finitely many values in, say, $\mathcal X_1 \times \cdots \times \mathcal X_d$,
\citet{LoWeSuGe13} consider an infinite ensemble of totally randomized and fully developed trees. At each cell, the ${\ell}$-th tree is grown by selecting a variable $X^{(j)}$ uniformly among the features that have not been used in the parent nodes, and by subsequently dividing the cell into $|\mathcal{X}_j|$ children (so the number of children equals the number of modalities of the selected variable). In this framework, it can be shown that the population version of MDI$(X^{(j)}$) computed with the whole forest satisfies
\begin{align*}
\mbox{MDI}^{\star}(X^{(j)}) = \sum_{k=0}^{p-1} \frac{1}{{k \choose p} (p-k)} \sum_{B \in \mathcal{P}_k(V^{-j})} I(X^{(j)}; Y|B),
\end{align*}
where $V^{-j}=\{1, \hdots, j-1, j+1, \hdots, p\}$, $\mathcal{P}_k(V^{-j})$  the set of subsets of $V^{-j}$ of cardinality $k$, and $I(X^{(j)}; Y |B)$  the {\it conditional mutual information} of $X^{(j)}$ and $Y$ given the variables in $B$. In addition,
$$
\sum_{j = 1}^p \mbox{MDI}^{\star}(X^{(j)}) = I(X^{(1)}, \hdots, X^{(p)}; Y).
$$
These results show that the information $I(X^{(1)}, \hdots, X^{(p)}; Y)$ is the sum of the importances of each variable, which can itself be made explicit using the information values $I(X^{(j)}; Y|B)$ between each variable $X^{(j)}$ and the output $Y$, conditional on variable subsets $B$ of different sizes.

\citet{LoWeSuGe13} define a variable $X^{(j)}$ as irrelevant with respect to $B \subset V=\mathcal X_1\times \cdots \times \mathcal X_p$ whenever $I(X^{(j)};Y|B)=0$. Thus, $X^{(j)}$ is irrelevant with respect to $V$ if and only if $\mbox{MDI}^{\star}(X^{(j)})=0$. It is easy to see that if an additional irrelevant variable $X^{(p+1)}$ is added to the list of variables, then, for any $j$, the variable importance $\mbox{MDI}^{\star}(X^{(j)})$ computed with a single tree does not change if the tree is built with the new collection of variables $V \cup \{X^{(p+1)}\}$. In other words, building a tree with an additional irrelevant variable does not alter the importances of the other variables in an infinite sample setting.

The most notable results regarding $\mbox{MDA}^{\star}$ are due to \citet{Is07}, who studies a slight modification of the criterion replacing permutation by feature noising. To add noise to a variable $X^{(j)}$, one considers a new observation $\bX$, take $\bX$ down the tree and stop when a split is made according to the variable $X^{(j)}$. Then the right or left child node is selected with probability $1/2$, and this procedure is repeated for each subsequent node (whether it is performed along the variable $X^{(j)}$ or not). The variable importance $\mbox{MDA}^{\star}(X^{(j)})$ is still computed by comparing the error of the forest with that of the ``noisy'' forest. Assuming that the forest is consistent and that the regression function is piecewise constant, \citet{Is07} gives the asymptotic behavior of $\mbox{MDA}^{\star}(X^{(j)})$ when the sample size tends to infinity. This behavior is intimately related to the set of subtrees (of the initial regression tree) whose roots are split along the coordinate $X^{(j)}$.

Let us lastly mention the approach of \citet{GrBeSa13}, who compute the MDA criterion for several distributions of $(\bX,Y)$. For example, consider a model of the form
\begin{align*}
Y = m(\bX) + \varepsilon,
\end{align*}
where $(\bX, \varepsilon)$ is a Gaussian random vector, and assume that the correlation matrix $C$ satisfies $C = [\mbox{Cov}(X^{(j)}, X^{(k)})]_{1\leq j,k\leq p} = (1-c) I_p + c \mathds{1} \mathds{1}^\top$ (the symbol $\top$ denotes transposition, $\mathds 1=(1, \hdots, 1)^{\top}$, and $c$ is a constant in $(0,1)$). Assume, in addition, that $\mbox{Cov}(X^{(j)},Y) = \tau_0$ for all $j \in \{1, \hdots, p\}$. Then, for all $j$,
$$
\mbox{MDA}^{\star}(X^{(j)}) = 2 \bigg( \frac{\tau_0}{1 - c + pc}\bigg)^2.
$$
Thus, in the Gaussian setting, the variable importance decreases as the inverse of the square of $p$ when the number of correlated variables $p$ increases.

\subsection{Related works}

The empirical properties of the \mbox{MDA} criterion have been extensively explored and compared in the statistical computing literature.  Indeed, \citet{ArKi08}, \citet{StBoKnAuZe08}, \citet{NiMa09}, \citet{AuAl11}, and \citet{ToLe11} stress the negative effect of correlated variables on MDA performance. In this respect, \citet{genuer} noticed that \mbox{MDA} is less able to detect the most relevant variables when
the number of correlated features increases. Similarly, the empirical study of \citet{ArKi08} points out that both \mbox{MDA} and \mbox{MDI} behave poorly when  correlation increases---these results have been experimentally confirmed by \citet{AuAl11} and \citet{ToLe11}.
An argument of \citet{StBoKnAuZe08} to justify the bias of \mbox{MDA} in the presence of correlated variables is that the algorithm evaluates the marginal importance of the variables instead of taking into account their effect conditional on each other. A way to circumvent this issue is to combine random forests and the \textit{Recursive Feature Elimination} algorithm of \citet{GuWeBaVa02}, as in \citet[][]{GrBeSa13}. Detecting relevant features can also be achieved via hypothesis testing \citep[][]{MeHo14a}---a principle that may be used to detect more complex structures of the regression function, like for instance its additivity \citep[][]{MeHo14}.

\section{Extensions}
{\bf Weighted forests.} In Breiman's \citeyearpar{Br01} forests, the final prediction is the average of the individual tree outcomes. A natural way to improve the method is to incorporate tree-level weights to emphasize more accurate trees in prediction \citep[][]{WiFrBi13}. A closely related idea, proposed by \citet{BeAdHe12}, is to guide  tree building---via resampling of the training set and other ad hoc randomization procedures---so that each tree will complement as much as possible the existing trees in the ensemble. The resulting {\it Dynamic Random Forest} (DRF) shows significant improvement in terms of accuracy on $20$ real-based data sets compared to the standard, static, algorithm.

{\bf Online forests.} In its original version, random forests is an {\it offline algorithm}, which is given the whole  data set from the beginning and  required to output an answer. In contrast, {\it online algorithms} do not require that the entire training set is accessible  at  once.
These  models  are  appropriate  for streaming  settings,  where  training  data  is  generated over time and must be incorporated into the model as quickly as possible. Random forests have been extended to the online framework in several ways \citep[][]{SaLeSaGoBi09,DeMaFr13,LaRoTe14}. In \citet{LaRoTe14},  so-called {\it Mondrian forests} are grown in an online fashion and achieve competitive predictive performance comparable with other online random forests while being faster. When building online forests, a major difficulty is to decide when the amount of data is sufficient to cut a cell. Exploring this idea, \citet{YiSoDeZh12} propose  {\it Information Forests}, whose construction consists in deferring classification until a measure of {\it classification confidence} is sufficiently high, and in fact break down the data so as to maximize this measure. An interesting theory related to these greedy trees can be found in \citet{BiDe13}.

{\bf Survival forests.} Survival analysis attempts to deal with analysis of time duration until one or more events happen. Most often, survival analysis is also concerned with incomplete data, and particularly right-censored data, in fields such as clinical trials. In this context, parametric approaches such as proportional hazards are commonly used, but fail to model nonlinear effects. Random forests have been extended to  the survival context by \citet{IsKoBlLa08}, who prove consistency of {\it Random Survival Forests} (RSF) algorithm assuming that all variables are categorical. \citet{YaWaFa10} showed that by incorporating kernel functions into RSF, their algorithm KIRSF achieves better results in many situations.
\citet{IsKoChMi11} review the use of the {\it minimal depth}, which measures the predictive quality of variables in survival trees.

{\bf Ranking forests.} \citet{ClDeVa13} have extended random forests to deal with ranking problems and propose an algorithm called \emph{Ranking Forests} based on the ranking trees of \citet{ClVa09}. Their approach relies on nonparametric scoring and ROC curve optimization in the sense of the AUC criterion.

{\bf Clustering forests. }\citet{YaChJo13} present a new clustering ensemble method called {\it Cluster Forests} (CF) in the context of unsupervised classification. CF randomly probes a high-dimensional data cloud to obtain good local clusterings,  then aggregates via spectral clustering to obtain cluster assignments for the whole data set. The search for good local clusterings is guided by a cluster quality measure, and CF progressively improves each local clustering in a fashion that resembles tree growth in random forests.

{\bf Quantile forests.}
\citet{Me06} shows that random forests provide information about the full conditional distribution of the response variable, and thus can be used for quantile estimation.

{\bf Missing data.} One of the strengths of random forests is that they can handle missing data. The procedure, explained in \citet{Br03}, takes advantage of the so-called {\it proximity matrix}, which measures the proximity between pairs of observations in the forest, to estimate  missing values. This measure is the empirical counterpart of the kernels defined in Section $3.2$. Data imputation based on random forests has further been explored by \citet{RiHoSt10}, \citet{CrFi08}, and extended to unsupervised classification by \citet{Is13Imp}.

{\bf Single class data.} One-class classification is a binary classification task for which only one class of samples is available for learning. \citet{DeBePeHe13} study the {\it One Class Random Forests} algorithm, which is designed to solve this particular problem. \citet{GeMeAy13} have introduced a supervised learning algorithm called {\it Spatially Adaptive Random Forests} to deal with semantic image segmentation applied to medical imaging protocols.
 Lastly, in the context of multi-label classification, \citet{JoGeWe14} adapt the idea of random projections applied to the output space to enhance tree-based ensemble methods by improving accuracy while  significantly reducing the computational burden.

{\bf Unbalanced data set.}
Random forests can naturally be adapted to fit the unbalanced data framework by down-sampling the majority class and growing each tree on a more balanced data set \citep[][]{BrChLi04, KuJo13}.
An interesting application in which unbalanced data sets are involved is by \citet{FiHoZuWiShMuHoRiShKe10}, who explore the continent-wide inter-annual migrations of common North American birds. They use random forests for which each tree is trained and allowed to predict on a particular (random) region in space and time.

\section{Conclusion and perspectives}

The authors trust that this review paper has provided an overview of some of the recent literature on random forests and offered insights into how new and emerging fields are impacting the method. As statistical applications become increasingly sophisticated, massive and complex data sets require today the development of algorithms that ensure global competitiveness, achieving both computational efficiency and safe with high-dimension models and huge number of samples. It is our belief that forests and their basic principles (``divide and conquer'', resampling, aggregation, random search of the feature space) offer simple but fundamental ideas that may leverage new state-of-the-art algorithms.

It remains however that the present results are insufficient to explain in full generality the remarkable behavior of random forests. The authors' intuition is that tree aggregation models are able to estimate patterns that are more complex than classical ones---patterns that cannot be simply characterized by standard sparsity or smoothness conditions. These patterns, which are beyond the reach of classical methods, are still to be discovered, quantified, and mathematically described.

It is sometimes alluded to that random forests have the flavor of deep network architectures \citep[e.g.,][]{Be09}, insofar as ensemble of trees allow to discriminate between a very large number of regions. Indeed, the identity of the leaf node with which a data point is associated for each tree forms a tuple that can represent a considerable quantity of possible patterns, because the total intersections of the leaf regions can be exponential in the number of trees. This point of view, largely unexamined, could be one of the reasons for the success of forests on large-scale data. As a matter of fact, the connection between random forests and neural networks is largely unexamined \citep[][]{We14}.

Another critical issue is how to choose tuning parameters that are optimal in a certain sense, especially the size $a_n$ of the preliminary resampling. By default, the algorithm runs in bootstrap mode (i.e., $a_n=n$ points selected with replacement) and although this seems to give excellent results, there is to date no theory to support this choice. Furthermore, although random forests are fully grown in most applications, the impact of tree depth on the statistical performance of the algorithm is still an open question.

\paragraph{Acknowledgments.} We thank the Editors and three anonymous referees for valuable comments and insightful suggestions.
\bibliography{biblio-TEST}
\end{document}